\input amstex
\documentstyle{amsppt}
%----------------------------------------------------------------
% Title:     Direct and inverse conversion formulas associated
%            with Khabibullin's conjecture for integral inequalities.
% Author:    Ruslan Sharipov
% Comments:  AmSTeX, 7 pages, amsppt style
% MSC-class: 26D10, 26D15, 39B62, 47A63
%----------------------------------------------------------------
%           Replacement for output macro definition
%
\catcode`@=11
\redefine\output@{%
  \def\break{\penalty-\@M}\let\par\endgraf
  \ifodd\pageno\global\hoffset=105pt\else\global\hoffset=8pt\fi  
  \shipout\vbox{%
    \ifplain@
      \let\makeheadline\relax \let\makefootline\relax
    \else
      \iffirstpage@ \global\firstpage@false
        \let\rightheadline\frheadline
        \let\leftheadline\flheadline
      \else
        \ifrunheads@ %\let\makefootline\relax
        \else \let\makeheadline\relax
        \fi
      \fi
    \fi
    \makeheadline \pagebody \makefootline}%
  \advancepageno \ifnum\outputpenalty>-\@MM\else\dosupereject\fi
}
\def\Beta{\mathchar"0\hexnumber@\rmfam 42}
\catcode`\@=\active
%----------------------------------------------------------------
\nopagenumbers
\def\negskp{\hskip -2pt}

%\def\tr{\operatorname{tr}}

%\accentedsymbol\bd{\kern 2pt\bar{\kern -2pt d}}
%\accentedsymbol\bbd{\kern 2pt\bar{\kern -2pt\bold d}}
%\def\vtrule{\vrule height 12pt depth 6pt}
%\def\vtttrule{\vrule height 12pt depth 19pt}
%\def\boxit#1#2{\vcenter{\hsize=122pt\offinterlineskip\hrule
%  \line{\vtttrule\hss\vtop{\hsize=120pt\centerline{#1}\vskip 5pt
%  \centerline{#2}}\hss\vtttrule}\hrule}}
\def\blue#1{#1}

\catcode`#=11\def\diez{#}\catcode`#=6
\catcode`_=11\def\podcherkivanie{_}\catcode`_=8
%\catcode`~=11\def\volna{~}\catcode`~=\active
\def\mycite#1{\cite{\blue{#1}}\immediate\special{ps:
     ShrHPSdict begin /ShrBORDERthickness 0 def}}

\def\mytag#1{%
    \tag#1}
\def\mythetag#1{\thetag{\blue{#1}}\immediate\special{ps:
     ShrHPSdict begin /ShrBORDERthickness 0 def}}
\def\myrefno#1{\no#1}
\def\myhref#1#2{\blue{#2}\immediate\special{ps:
     ShrHPSdict begin /ShrBORDERthickness 0 def}}
\def\myEarXivlink{\myhref{http://arXiv.org}{http:/\negskp/arXiv.org}}

\def\mytheorem#1{\csname proclaim\endcsname{Theorem #1}}
\def\mytheoremwithtitle#1#2{\csname proclaim\endcsname{Theorem #1#2}}
\def\mythetheorem#1{\blue{#1}\immediate\special{ps:
     ShrHPSdict begin /ShrBORDERthickness 0 def}}
\def\mylemma#1{\csname proclaim\endcsname{Lemma #1}}
\def\mylemmawithtitle#1#2{\csname proclaim\endcsname{Lemma #1#2}}
\def\mythelemma#1{\blue{#1}\immediate\special{ps:
     ShrHPSdict begin /ShrBORDERthickness 0 def}}
\def\mycorollary#1{\csname proclaim\endcsname{Corollary #1}}

\def\myconjecture#1{\csname proclaim\endcsname{Conjecture #1}}
\def\myconjecturewithtitle#1#2{\csname proclaim\endcsname{Conjecture #1#2}}
\def\mytheconjecture#1{\blue{#1}\immediate\special{ps:
     ShrHPSdict begin /ShrBORDERthickness 0 def}}

%----------------------------------------------------------------
% Cyrillic fonts definition
%\font\eightcyr=wncyr8
%----------------------------------------------------------------
\pagewidth{360pt}
\pageheight{606pt}
\topmatter
\title
Direct and inverse conversion formulas associated
with Khabibullin's conjecture for integral inequalities.
\endtitle
\rightheadtext{Direct and inverse conversion formulas \dots}
\author
R.~A.~Sharipov
\endauthor
\address 5 Rabochaya street, 450003 Ufa, Russia\newline
\vphantom{a}\kern 12pt Cell Phone: +7(917)476 93 48
\endaddress
\email \vtop to 30pt{\hsize=280pt\noindent
\myhref{mailto:r-sharipov\@mail.ru}
{r-sharipov\@mail.ru}\newline
\myhref{mailto:R\podcherkivanie Sharipov\@ic.bashedu.ru}
{R\_\hskip 1pt Sharipov\@ic.bashedu.ru}\vss}
\endemail
\urladdr
\vtop to 20pt{\hsize=280pt\noindent
\myhref{http://ruslan-sharipov.ucoz.com/}
{http:/\negskp/ruslan-sharipov.ucoz.com}\newline
\myhref{http://www.freetextbooks.narod.ru/}
{http:/\negskp/www.freetextbooks.narod.ru}\newline
\myhref{http://sovlit2.narod.ru/}
{http:/\negskp/sovlit2.narod.ru}\vss}
\endurladdr
\abstract
    Khabibullin's conjecture deals with two linear integral inequalities
for some non-negative continuous function $q(t)$. The integral in the first
of these two inequalities converts $q(t)$ into another function of one 
variable $g(t)$. This integral yields the direct conversion formula. An
inverse conversion formula means a formula expressing $q(t)$ back through
$g(t)$. Such an inverse conversion formula is derived. 
\endabstract
\endtopmatter
%\loadbold
%\loadeufb
\TagsOnRight
\document

\head
1. Introduction.
\endhead
     In this paper the following statement of Khabibullin's conjecture is 
used. \myconjecturewithtitle{1.1}{ (Khabibullin)} Let $\alpha>0$ and let 
$q=q(t)$ be a continuous function such that $q(t)\geqslant 0$ for all $t>0$. 
Then the inequality 
$$
\hskip -2em
\int\limits^{\,1}_0\left(\,\,\int\limits^{\,1}_x(1-y)^{n-1}
\,\frac{dy}{y}\right)q(tx)\,dx\leqslant t^{\alpha-1}
\mytag{1.1}
$$
fulfilled for all \kern 2pt $0\leqslant t<+\infty$ implies the inequality
$$
\int\limits^{+\infty}_0 q(t)\,\ln\Bigl(1+\frac{1}{t^{\,2\kern 0.2pt\alpha}}
\Bigr)\,dt\leqslant\pi\,\alpha\prod^{n-1}_{k=1}\Bigl(1+\frac{\alpha}{k}\Bigr).
$$
\endproclaim
     Initially,  the conjecture~\mytheconjecture{1.1} was formulated in 
\mycite{1} and \mycite{2}, though in some different form. Later in 
\mycite{3} it was reformulated in a form very close to the above statement. 
In \mycite{4} the conjecture~\mytheconjecture{1.1} was proved to be valid 
for $0<\alpha\leqslant 1/2$. Another proof of this result was given
in \mycite{5}.\par
      The approach of the paper \mycite{5} is based on the kernel function 
$A_n(x)$. The kernel function $A_n(x)$ is defined by the inner integral in 
the formula \mythetag{1.1}:
$$
\hskip -2em
A_n(x)=\int\limits^{\,1}_{\!x}(1-y)^n\,\frac{dy}{y}.
\mytag{1.2}
$$
In terms of the kernel function \mythetag{1.2} the inequality \mythetag{1.1} 
is written as 
$$
\hskip -2em
\int\limits^{\,1}_0\!A_{n-1}(x)\,q(t\,x)\,dx\leqslant t^{\alpha-1}\text{, 
\ where \ }t\geqslant 0.
\mytag{1.3}
$$
By changing the variable $x$ for the variable $y=t\,x$ in the integral 
\mythetag{1.3} we get
$$
\hskip -2em
\int\limits^{\,t}_0\!A_{n-1}(y/t)\,q(y)\,dy\leqslant t^\alpha\text{, 
\ where \ }t>0.
\mytag{1.4}
$$
The value $t=0$ is an exception when transforming \mythetag{1.3} into 
\mythetag{1.4}. We omit this exceptional value from our further 
considerations.\par
      Looking at the left hand side of the inequality \mythetag{1.4}, 
we define the following integral transformation that converts a function 
$q=q(t)$ defined on the half-line $t>0$ into another function $g=g(t)$ 
defined on the same half-line $t>0$:
$$
\hskip -2em
g(t)=\int\limits^{\,t}_0\!A_n(y/t)\,q(y)\,dy.
\mytag{1.5}
$$
The formula \mythetag{1.5} is called the direct conversion formula. The 
main goal of this paper is to derive an inverse conversion formula that 
converts $g(t)$ back to $q(t)$. 
\head
2. Properties of the kernel function. 
\endhead
\parshape 16 180pt 180pt 180pt 180pt 180pt 180pt 180pt 180pt 
180pt 180pt 180pt 180pt 180pt 180pt 180pt 180pt 180pt 180pt 
180pt 180pt 180pt 180pt 180pt 180pt 180pt 180pt 180pt 180pt 
180pt 180pt 0pt 360pt 
     The kernel function $A_n(x)$ used in the direct conversion formula 
\mythetag{1.5}. Its properties were studied in \mycite{5}. 
\vadjust{\vskip 5pt\hbox to 0pt{\kern 
0pt \includegraphics{khab02_01.eps}\hss}\vskip -5pt}This is 
a decreasing smooth function on the segment $(0,\,1]$ vanishing at the 
point $x=1$ and having the logarithmic singularity 
$$
\hskip -2em
A_n(x)\sim -\ln x
\mytag{2.1}
$$
at the point $x=0$. Its graph is shown on Fig\.~2.1. There are two 
explicit formulas for $A_n(x)$. Here is the first of them
$$
\hskip -2em
A_n(x)=\kern -4pt\sum^\infty_{m=n+1}\kern -4pt\frac{(1-x)^m}{m}.
\mytag{2.2}
$$
From \mythetag{2.2} we immediately derive the following vanishing 
conditions:
$$
\hskip -2em
\frac{d^{\kern 0.2pt k}\!A_n(x)}{dx^k}\,{\vrule height 16pt depth 7pt}_{\ x=1}
\!=0\text{\ \ for \ }k=0,\,1,\,\ldots,\,n.
\mytag{2.3}
$$
The sum in the second explicit formula for the kernel function $A_n(x)$ is 
finite
$$
\hskip -2em
A_n(x)=-\ln x-\sum^n_{m=1}\frac{(1-x)^m}{m}.
\mytag{2.4}
$$
The formula \mythetag{2.1} is immediate from \mythetag{2.4}.\par
     Note that the function $A_n(x)$ enters the formula \mythetag{1.5} in 
the form of $A_n(y/t)$ with the composite argument $x=y/t$. Substituting 
$x=y/t$ into \mythetag{2.4}, we get 
$$
\hskip -2em
A_n(y/t)=-\ln y+\ln t-\sum^n_{m=1}\frac{(1-y/t)^m}{m}.
\mytag{2.5}
$$
The function \mythetag{2.5} can be treated as a function of two variables 
$y$ and $t$. The partial derivatives of $A_n(y/t)$ with respect to $y$ and 
$t$ can be calculated explicitly:
$$
\xalignat 2
&\hskip -2em
\frac{\partial A_n(y/t)}{\partial y}=-\frac{(t-y)^n}{t^n\,y},
&&\frac{\partial A_n(y/t)}{\partial t}=\frac{(t-y)^n}{t^{n+1}}.
\mytag{2.6}
\endxalignat
$$
The formulas \mythetag{2.6} are easily derived from the following 
formula:
$$
\hskip -2em
\frac{dA_n(x)}{dx}=\frac{-(1-x)^n}{x}.
\mytag{2.7}
$$
As for the formula \mythetag{2.7}, it is derived from \mythetag{2.4} by 
means of direct calculations. The reader can find more details in \mycite{5}. 
\head
3. The first derivative of the function $g(t)$. 
\endhead
      The integral in the direct conversion formula \mythetag{1.5} is 
assumed to be finite for at least one value of $t=t_0>0$ as a convergent 
improper Riemann integral or as a Lebesque integral. Under this assumption 
we have the following lemma.
\mylemma{3.1} Let $q=q(t)$ be a non-negative continuous function on the 
open half-line $t>0$, i\.\,e\. $q(t)\geqslant 0$ for all $t>0$. If the 
integral \mythetag{1.5} is finite for some $t_0>0$ then for all $t>0$ the 
following integrals are finite:
$$
\xalignat 2
&\hskip -2em
\int\limits^{\,t}_0 q(y)\,dy<\infty,
&&\int\limits^{\,t}_0 |\ln y|\,q(y)\,dy<\infty.
\mytag{3.1}
\endxalignat
$$
\endproclaim
\demo{Proof} Since $q=q(t)$ is a continuous function, in order to prove 
the inequalities \mythetag{3.1} for all $t>0$ it is sufficient to prove 
them for some particular $t=y_0>0$. 
Since
$$
\hskip -2em
\int\limits^{\ t_0}_0\!A_n(y/t_0)\,q(y)\,dy<\infty,
\mytag{3.2}
$$
we choose $t=t_0$ and from \mythetag{2.4} we derive $A_n(y/t_0)\to +\infty$ 
as $y\to +0$. Hence, there is some $y_0>0$ such that $A_n(y/t_0)>1$ for all 
$0<y\leqslant y_0$. Multiplying by $q(y)$ and taking into account that 
$q(y)\geqslant 0$, from $A_n(y/t_0)>1$ we derive
$$
\hskip -2em
q(y)\leqslant A_n(y/t_0)\,q(y)\text{\ \ for \ }0<y\leqslant y_0.
\mytag{3.3}
$$
Integrating the inequality \mythetag{3.3} we get
$$
\hskip -2em
\int\limits^{\ y_0}_0 q(y)\,dy\leqslant\int\limits^{\ y_0}_0\!A_n(y/t_0)
\,q(y)\,dy.
\mytag{3.4}
$$
Note that the integration interval in \mythetag{3.4} differs from that 
of \mythetag{3.2}. However, since both $A_n(y/t_0)$ and $q(y)$ are continuous 
functions, extending or shrinking the integration interval does not affect 
the finiteness of the integral \mythetag{3.2}. Therefore, combining 
\mythetag{3.2} and \mythetag{3.4}, we get the inequality
$$
\int\limits^{\ y_0}_0 q(y)\,dy<\infty.
$$
Thus, the first inequality \mythetag{3.1} of the lemma~\mythelemma{3.1} is 
proved.\par
     In order to prove the second inequality we use the formula 
\mythetag{2.4} again. Applying this formula to the function $A_n(y/t_0)$, 
we obtain
$$
\hskip -2em
\lim_{y\to +0}\frac{A_n(y/t_0)}{|\ln y|}=1. 
\mytag{3.5}
$$ 
The equality \mythetag{3.5} means that there is some $y_0>0$ such that 
$$
\hskip -2em
|\ln y|<2\,A_n(y/t_0)\text{\ \ for all \ }0<y\leqslant y_0.
\mytag{3.6}
$$
Multiplying \mythetag{3.6} by $q(y)$ and taking into account that 
$q(y)\geqslant 0$, we get
$$
\hskip -2em
|\ln y|\,q(y)\leqslant 2\,A_n(y/t_0)\,q(y)\text{\ \ for all \ }
0<y\leqslant y_0.
\mytag{3.7}
$$
Now, integrating the inequality \mythetag{3.7}, we obtain
$$
\hskip -2em
\int\limits^{\ y_0}_0|\ln y|\,q(y)\,dy\leqslant 2\!\int\limits^{\ y_0}_0
\!A_n(y/t_0)\,q(y)\,dy.
\mytag{3.8}
$$
Combining \mythetag{3.2} and \mythetag{3.8}, we derive
$$
\hskip -2em
\int\limits^{\ y_0}_0|\ln y|\,q(y)\,dy<\infty. 
\mytag{3.9}
$$
The second inequality \mythetag{3.1} of the lemma~\mythelemma{3.1} is also
proved. As we already said above, the difference in upper limits of the 
integrals \mythetag{3.1} and \mythetag{3.9} does not matter for finiteness 
of the integral \mythetag{3.9}.
\qed\enddemo
     Let's return back to the direct conversion formula \mythetag{1.5} and,
choosing some constant $b$, let's subdivide the integral \mythetag{1.5} 
into two integrals: 
$$
\xalignat 2
&\hskip -2em
I_1(t)=\int\limits^{\,b}_0\!A_n(y/t)\,q(y)\,dy,
&&I_2(t)=\int\limits^{\,t}_b\!A_n(y/t)\,q(y)\,dy.
\quad
\mytag{3.10}
\endxalignat
$$
Writing \mythetag{3.10}, we assume that $0<b<t$. Let $c$ be another constant 
such that $0<b<t<c$. The first integral \mythetag{3.10} is an improper 
integral with the singularity at its lower limit $y=0$. The second integral 
\mythetag{3.10} is a proper integral. In addition to \mythetag{3.10}, let's 
consider the following two integrals:
$$
\xalignat 2
&\hskip -2em
I_3(t)=\int\limits^{\,b}_0\frac{\partial A_n(y/t)}{\partial t}\,q(y)\,dy,
&&I_4(t)=\int\limits^{\,t}_b\frac{\partial A_n(y/t)}{\partial t}\,q(y)\,dy.
\quad
\mytag{3.11}
\endxalignat
$$
Due to \mythetag{2.5} and \mythetag{2.6} the functions 
$$
\xalignat 2
&\hskip -2em
A_n(y/t)\,q(y),
&&\frac{\partial A_n(y/t)}{\partial t}\,q(y)
\mytag{3.12}
\endxalignat
$$
both are functions of two variables $y$ and $t$ which are continuous within 
the closed rectangle $R_2=\{(y,t)\in\Bbb R^2\!:\,b\leqslant y\leqslant c,\ 
b\leqslant t\leqslant c\}$. Therefore we can apply the theorem~4' from 
\S\,53 of Chapter~\uppercase\expandafter{\romannumeral 6} in \mycite{6} to 
the integrals $I_2(t)$ and $I_4(t)$. This theorem says that $I_2(t)$ is a 
differentiable function such that
$$
\hskip -2em
\frac{dI_2(t)}{dt}=I_4(t)+A_n(t/t)\,q(t).
\mytag{3.13}
$$
Note that $t/t=1$ and, according to \mythetag{2.3}, $A_n(1)=0$. Therefore 
\mythetag{3.13} reduces to
$$
\hskip -2em
\frac{dI_2(t)}{dt}=I_4(t).
\mytag{3.14}
$$\par
     Now let's proceed to the integrals $I_1(t)$ and $I_3(t)$ in 
\mythetag{3.10} and \mythetag{3.11}. The functions \mythetag{3.12} both 
are continuous functions of two variables within the semi-open rectangle 
$R_1=\{(y,t)\in\Bbb R^2\!:\,0<y\leqslant b,\ b\leqslant t\leqslant c\}$. 
Using \mythetag{2.5} and \mythetag{2.6}, one can easily prove that there 
are two constants $C_1$ and $C_2$ such that
$$
\xalignat 2
&\hskip -2em
|A_n(y/t)\,q(y)|\leqslant C_1\,(|\ln y|+1)\,q(y),
&&\left|\frac{\partial A_n(y/t)}{\partial t}\,q(y)\right|\leqslant C_2
\,q(y)
\qquad
\mytag{3.15}
\endxalignat
$$
for all $(y,t)$ within the semi-open rectangle $R_1$. Due to \mythetag{3.15} 
and \mythetag{3.1} we can apply the theorem~1 from \S\,54 of 
Chapter~\uppercase\expandafter{\romannumeral 6} in \mycite{6} to the improper
integrals $I_1(t)$ and $I_3(t)$ in \mythetag{3.10} and \mythetag{3.11}. This 
theorem says that both of these two improper integrals converge uniformly in 
$t$ over the interval $b\leqslant t\leqslant c$. Due to the uniform 
convergence, we can apply the theorem~8 from \S\,54 of 
Chapter~\uppercase\expandafter{\romannumeral 6} in \mycite{6} to $I_1(t)$ 
and $I_3(t)$. This theorem says that $I_1(t)$ is a differentiable function 
and 
$$
\hskip -2em
\frac{dI_1(t)}{dt}=I_3(t).
\mytag{3.16}
$$\par
     Since $g(t)=I_1(t)+I_2(t)$, now we can combine the formulas 
\mythetag{3.14} and \mythetag{3.16} and derive the following formula for 
the first derivative of $g(t)$:
$$
\pagebreak
\hskip -2em
g'(t)=\int\limits^{\,t}_0\frac{\partial A_n(y/t)}{\partial t}\,dy.
\mytag{3.17}
$$
More precisely this result is formulated in the following theorem.
\mytheorem{3.1} If $q=q(t)$ is a non-negative continuous function on the 
half-line $t>0$ and if the integral \mythetag{1.5} is finite for at least 
one value $t=t_0>0$, then $g(t)$ is a differentiable function and its 
derivative is given by the formula \mythetag{3.17}.
\endproclaim
\head
4. Higher order derivatives of the function $g(t)$. 
\endhead
     Let's apply the second formula \mythetag{2.6} to \mythetag{3.17}. Then 
the formula \mythetag{3.17} for $g'(t)$ is brought to the following more 
explicit form:
$$
\hskip -2em
t^{n+1}\,g'(t)=\int\limits^{\,t}_0(t-y)^n\,q(y)\,dy. 
\mytag{4.1}
$$
Let's denote $\tilde g(t)=t^{n+1}\,g'(t)$ and write \mythetag{4.1} as
$$
\hskip -2em
\tilde g(t)=\int\limits^{\,t}_0(t-y)^n\,q(y)\,dy. 
\mytag{4.2}
$$
The formula is quite similar to \mythetag{1.5}. It is even simpler than 
\mythetag{1.5} since the function $(t-y)^n$ has no logarithmic singularity, 
though this fact does not matter for us. Applying the same arguments as 
in proving the theorem~\mythetheorem{3.1} above, we can prove the following 
theorem for $\tilde g(t)$. 
\mytheorem{4.1} If $q=q(t)$ is a non-negative continuous function on the 
half-line $t>0$ and if the integral \mythetag{1.5} is finite for at least 
one value $t=t_0>0$, then $\tilde g(t)$ in \mythetag{4.2} is a differentiable 
function and its derivative is given by the formula
$$
\hskip -2em
\tilde g'(t)=\int\limits^{\,t}_0n\,(t-y)^{n-1}\,q(y)\,dy. 
\mytag{4.3}
$$
\endproclaim
The next theorem yields the second derivative of the function $\tilde g(t)$.
\mytheorem{4.2} If $q=q(t)$ is a non-negative continuous function on the 
half-line $t>0$ and if the integral \mythetag{1.5} is finite for at least 
one value $t=t_0>0$, then the function $\tilde g(t)$ in \mythetag{4.2} is 
a twice differentiable function and its second order derivative is given 
by the formula
$$
\hskip -2em
\tilde g''(t)=\int\limits^{\,t}_0n\,(n-1)\,(t-y)^{n-2}\,q(y)\,dy. 
\mytag{4.4}
$$
\endproclaim
Acting repeatedly, one can prove the series of theorems saying that 
$\tilde g(t)$ is an $n$ times differentiable function and derive the 
formula generalizing \mythetag{4.3} and \mythetag{4.4}:
$$
\pagebreak
\hskip -2em
\frac{d^{\kern 0.2pt k}\tilde g(t)}{dt^k}=\frac{n!}{(n-k)!}
\int\limits^{\,t}_0\!(t-y)^{n-k}\,q(y)\,dy,\qquad k=0,\,\ldots,\,n. 
\mytag{4.5}
$$
For $k=n$ the formula \mythetag{4.5} reduces to the following one:
$$
\hskip -2em
\frac{d^{\kern 0.2pt n}\tilde g(t)}{dt^n}=n!\int\limits^{\,t}_0\!
q(y)\,dy. 
\mytag{4.6}
$$
Due to \mythetag{4.6} the function $g(t)$ in \mythetag{1.5} is an $(n+1)$ 
times differentiable function. 
\head
5. An inverse conversion formula. 
\endhead
     An inverse conversion formula is almost immediate from \mythetag{4.6}.
Indeed, due to \mythetag{3.1} the right hand side of the formula 
\mythetag{4.6} is a differentiable function. Differentiating \mythetag{4.6} 
and substituting $\tilde g(t)=t^{n+1}\,g'(t)$, we get
$$
\hskip -2em
q(t)=\frac{d^{\kern 0.3pt n+1}}{dt^{n+1}}\!\left(\frac{t^{n+1}
\,g'(t)}{n!}\right)\!.
\mytag{5.1}
$$
This is a required inverse conversion formula. Tne following ultimate theorem 
is associated with the formula \mythetag{5.1}.
\mytheorem{5.1} If $q=q(t)$ is a non-negative continuous function on the 
half-line $t>0$ and if the integral \mythetag{1.5} is finite for at least 
one value $t=t_0>0$, then the function $g(t)$ in \mythetag{1.5} is an 
$(n+2)$ times differentiable function such that $q(t)$ is expressed through 
its derivatives according to the formula \mythetag{5.1}.
\endproclaim
     The theorem~\mythetheorem{5.1} does not claim the formula \mythetag{5.1} 
to be the only way for expressing $q(t)$ through $g(t)$. But hopefully, the 
formula \mythetag{5.1} could be a useful tool for studying Khabibullin's 
conjecture~\mytheconjecture{1.1}.


\Refs
\ref\myrefno{1}\by Khabibullin~B.~N.\paper Paley problem for plurisubharmonic 
functions of a finite lower order\jour Mat\. Sbornik\vol 190\issue 2\yr 1999
\pages 145-157
\endref
\ref\myrefno{2}\by Khabibullin~B.~N.\paper The representation of a meromorphic 
function as a quotient of entire functions and the Paley problem in $\Bbb C^n$: survey of some results
\jour Mathematical Physics, Analysis, and geometry (Ukraine) \yr 2002\vol 9
\issue 2\pages 146-167\moreref see also
\myhref{http://arxiv.org/abs/math.CV/0502433}{math.CV/0502433} in Electronic 
Archive \myEarXivlink
\endref
\ref\myrefno{3}\by Khabibullin~B.~N.\paper A conjecture on some estimates 
for integrals\publ e-print \myhref{http://arXiv.org/abs/1005.3913/}
{arXiv:1005.3913} in Electronic Archive \myEarXivlink
\endref
\ref\myrefno{4}\by Baladai~R.~A, Khabibullin~B.~N.\paper Three equivalent 
conjectures on an estimate of integrals\publ e-print 
\myhref{http://arXiv.org/abs/1006.5140/}{arXiv:1006.5140} in Electronic 
Archive \myEarXivlink
\endref
\ref\myrefno{5}\by Sharipov~R.~A.\paper A note on Khabibullin's conjecture 
for integral inequalities\publ e-print 
\myhref{http://arXiv.org/abs/1008.0376/}{arXiv:1008} 
\myhref{http://arXiv.org/abs/1008.0376/}{.0376} 
in Electronic Archive \myEarXivlink
\endref
\ref\myrefno{6}\by Kudryavtsev~L.~D.\book Course of mathematical analysis. 
\rm Vol\.~\uppercase\expandafter{\romannumeral 2}\publ Visshaya Shkola 
publishers\publaddr Moscow\yr 1981
\endref
\endRefs
\enddocument
\end